\documentclass{amsart}
\usepackage{amsmath}
\usepackage{amsthm}
\usepackage{amsfonts}
\usepackage{amssymb}
\usepackage[all]{xy}
\usepackage{alltt}
\usepackage{multicol}
\usepackage{pdftexcmds}
\usepackage{etoolbox}
\usepackage{ifthen}
\usepackage[
    colorlinks=true,urlcolor=blue,linkcolor=blue,citecolor=blue
]{hyperref}
\usepackage{textcomp}
\usepackage{fonttable}
\usepackage{enumitem}
\usepackage{graphicx}
\usepackage{comment}
\usepackage{tikz}           
\usepackage{libertine}
\usepackage[libertine]{newtxmath}
\usepackage{centernot}
\usepackage[margin=1.5in]{geometry}

\theoremstyle{plain}
\newtheorem{theorem}{Theorem}

\newtheorem{conjecture}[theorem]{Conjecture}

\theoremstyle{definition}

\theoremstyle{definition}

\DeclareMathOperator{\rank}{rank}

\newcommand{\Z}{\ensuremath{\mathbf{Z}}}
\newcommand{\F}{\ensuremath{\ff}}
\newcommand{\R}{\ensuremath{\mathbf{R}}}
\newcommand{\Q}{\ensuremath{\mathbf{Q}}}

\newcommand{\C}{\ensuremath{\mathbf{C}}}

\newcommand{\n}{\noindent}
\newcommand{\map}[3]{\ensuremath{#1 : #2 \longrightarrow #3}}

\renewcommand{\F}{\ensuremath{\mathbf{F}}}

\begin{document}
 

\title{Is there an infinite field whose multiplicative group is indecomposable?}
\date{\today}
 
\author{Sunil K. Chebolu}
\address{Department of Mathematics \\
Illinois State University \\
Normal, IL 61790, USA}
\email{schebol@ilstu.edu}

\author{Keir Lockridge} 
\address {Department of Mathematics \\
Gettysburg College \\
Gettysburg, PA 17325, USA}
\email{klockrid@gettysburg.edu}


\thanks{The first author is supported by Simons Foundation: Collaboration Grant for Mathematicians (516354). }


\maketitle

\begin{abstract}
In \cite{cl-2}, we determined the finite fields with indecomposable multiplicative groups and conjectured that there is no infinite field whose multiplicative group is indecomposable. In this paper,  we prove this conjecture for several popular classes of fields including finitely generated fields, discrete valued fields, fields of Hahn series, local fields, global fields, and function fields. 
\end{abstract}

\section{introduction}

The purpose of this note is to address an interesting open question: {\em is there an infinite field whose multiplicative group is indecomposable?} A negative answer to this question would solve a 50-year-old unsolved problem of Fuchs, discussed below, concerning the realizability of indecomposable abelian groups over the class of fields!

 It is immediate from the definition of a field that its nonzero elements form an abelian group under multiplication, and these groups have been studied for decades. An abelian group is {\em indecomposable} if it is not isomorphic to a direct sum of two nontrivial groups. In the category of finitely generated abelian groups, the indecomposable groups are quite familiar (the cyclic groups of prime power or infinite order), but in the full category of abelian groups, they can be quite complicated. It is somewhat surprising that, as far as we know, our eponymous question has not been mentioned in the literature (other than our paper \cite{cl-2}), though it is an easy-to-pose question about fundamental mathematical objects.

The broader context for this problem is the study of units in rings. For a unital ring $R$, let $R^\times$ denote the group of units in $R$. The functor \[ \map{\mathcal{U}}{\mathbf{Rings}}{\mathbf{Groups}} \] determined by the assignment $\mathcal{U}(R) = R^\times$ is a tool for studying a facet of the interplay between groups and rings. A group in the image of \,$\mathcal{U}$ is called {\em realizable}. More than 50 years ago in \cite{fuchsprob}, L\'{a}szl\'{o} Fuchs posed the problem of determining the image of $\mathcal{U}$ when restricted to the class of commutative rings or fields, and many authors have studied the image or preimage of $\mathcal{U}$ when restricted to special classes of groups or rings. Such work usually falls into one of the following three categories; we give just a few examples in each case.
\begin{enumerate}
\item Given a special class of rings, compute the unit groups of its members. Examples: rings of integers in a number field (Dirichlet's Unit Theorem), finite rings (\cite{Dolzan}), group rings (\cite{higman}), torsion-free rings (\cite{dd}).
\item Given a special class of groups, determine the realizable members in this class. Examples: indecomposable abelian groups (\cite{cl-3}), dihedral groups (\cite{cl-4}), symmetric, alternating, and simple groups (\cite{do1}, \cite{do2}).
\item Given a special class of groups, determine the rings $R$ for which $R^\times$ is a member of this class. Example: cyclic groups (\cite{gilmer}, \cite{PS}).
\end{enumerate}

A pertinent example in the second category above is the problem of determining the indecomposable abelian groups that are realizable. We solved this problem in \cite{cl-3}; the next theorem gives a concise summary of the solution. Recall that a {\em Fermat prime} is a prime of the form $2^n + 1$ and a {\em Mersenne prime} is a prime of the form $2^n - 1$. A group is {\em torsion-free} if it contains no nontrivial elements of finite order.

\begin{theorem}[\cite{cl-3}]
Let $R$ be a unital ring whose group of units is an indecomposable abelian group.\label{ringind}
\begin{enumerate}
\item If $R^\times$ is finite, then $R^\times$ is isomorphic to $\C_1$, $\C_8$, $\C_{q-1}$ for $q$ a Fermat prime, or $\C_p$ for $p$ a Mersenne prime.
\item If $R^\times$ is infinite, then $R$ has characteristic 2 and $R^\times$ is an indecomposable torsion-free abelian group.
\end{enumerate}
Further, the groups $\C_1$, $\C_8$, $\C_{q-1}$ for $q$ a Fermat prime, and $\C_p$ for $p$ a Mersenne prime are each the group of units in some ring, and every indecomposable torsion-free abelian group is the group of units in a ring of characteristic 2.
\end{theorem}

Now suppose we wish to restrict our attention from rings to fields. One might expect this task to be easier or for the answer to fall out of the above theorem or its proof. For finite fields, this is indeed the case, and in \cite{cl-2} we proved the following theorem. Recall that, up to isomorphism, there is a unique field $\F_{p^n}$ of prime power order $p^n$, and every finite field has prime power order; further, the multiplicative group of a finite field is cyclic: $\F_{p^n}^\times \cong \C_{p^n - 1}$.

\begin{theorem}[{\cite[2.4]{cl-2}}] Let $F$ be a finite field. The multiplicative group $F^\times$ is indecomposable if and only if $F = \F_2, \F_9, \F_q$ where $q$ is a Fermat prime, or $\F_{p+1}$ where $p$ is a Mersenne prime. \end{theorem}

\n Note that the fields above account for all of the finite groups appearing in Theorem \ref{ringind}. Thus, it remains to determine which infinite fields have indecomposable multiplicative groups. We conjectured in \cite{cl-2} that no such field exists.

\begin{conjecture}[{\cite[4.5]{cl-2}}] \underline{There is no infinite field whose multiplicative group is indecomposable. }\label{mainc} \end{conjecture}

 We now prove this conjecture for several popular classes of fields including finitely generated fields, discrete valued fields, fields of Hahn series, local fields, global fields, and function fields. 
\n If true, Conjecture \ref{mainc} would complete the classification of fields with indecomposable multiplicative groups. It would also imply that the complete list of the indecomposable abelian groups that are realizable as the multiplicative group of a field is: $\C_1$, $\C_8$, $\C_{q-1}$ for $q$ a Fermat prime, and $\C_p$ for $p$ a Mersenne prime. 

In \cite{cl-2}, we proved that the multiplicative group of an infinite field whose characteristic is not 2 must be decomposable, so Conjecture \ref{mainc} is equivalent to the claim that the multiplicative group of an infinite field of characteristic 2 must be decomposable. If Conjecture \ref{mainc} is false, then what other features might a counterexample have? Are there important classes of fields we can quickly eliminate from consideration? Before investigating these questions, we will first summon a few basic facts from field theory in the next section; see \cite{Roman} for more details.

\section{Background from field theory}

A field might be built from its prime subfield by adjoining only algebraic elements, such as $\Q(\sqrt{2})$ (these are called {\em absolutely algebraic fields}). On the other hand, a field might be built from its prime subfield by adjoining only algebraically independent elements (they satisfy no polynomial relation), such as $\Q(\pi)$; such fields are examples of {\em purely transcendental extensions}. Other fields are obtained by combining algebraic and transcendental extensions. We would like to know how a possible counterexample to Conjecture \ref{mainc} might fit into this continuum of possibilities.

The purely transcendental extensions mentioned above are isomorphic to fields of rational functions, which we now define. Given a nonempty set $S$ of variables, and a field $B$, the {\em field of rational functions} $B(S)$ is the field of fractions for the ring $B[S]$ of polynomials with indeterminates in $S$ and coefficients in $B$; i.e., the elements of $B(S)$ are ratios $f/g$ of polynomials $f, g \in B[S]$ with $g \neq 0$. Since $B[x]$ is a unique factorization domain, we have that \[ B(S)^\times \cong B^\times \times \left(\bigoplus_{f \in \Delta} f^\Z\right) \cong B^\times \times \left(\bigoplus_{f \in \Delta} \Z\right),\] where $\Delta$ is the set of monic irreducible polynomials in $B[S]$ (see \cite{DummitFoote} for more details) and $f^\Z$ denotes the multiplicative infinite cyclic group generated by $f$. The {\em rank} of an abelian group $G$, written $\rank G$, is the size of a maximal linearly independent (over $\Z$) subset of $G$; from this definition, it is clear that $\rank H \leq \rank G$ whenever $H$ is a subgroup of $G$. Now, for any field $B$ and nonempty set $S$, we see that $B(S)^\times$ is a decomposable abelian group of infinite rank. 

Fields of rational functions are important in the current context because every field extension $A\, / \, B$ has a (possibly empty) transcendence basis. There is a subset $S\subseteq A$ of algebraically independent elements such that $A$ is an algebraic extension of the subfield $B(S)$ generated by $S$, which is field isomorphic to the field of rational functions over $B$ with indeterminates in $S$. For example, the extension $\Q(\pi)\,/\,\Q$ has transcendence basis $\{\pi\}$ and $\Q(\pi)$ is isomorphic to the field of rational functions $\Q(x)$.

We may now proceed with an analysis of Conjecture \ref{mainc} within this framework. 

\section{General features of a possible counterexample}

Suppose $F\, / \, \F_2$ is an infinite field whose multiplicative group is indecomposable. We immediately obtain that $F^\times$ is torsion-free by Theorem \ref{ringind}. This implies that $F$ contains no finite subfields other than $\F_2$, so $F$ contains no elements which are algebraic over $\F_2$ (put another way, $\F_2$ is algebraically closed in $F$). This eliminates, for example, the class of absolutely algebraic fields mentioned above. Now let $T$ be a transcendence basis for $F \, / \, \F_2$. The set $T$ is nonempty, for otherwise $F$ would be algebraic over $\F_2$. Since $F^\times \supseteq \F_2(T)^\times$ and $\F_2(T)^\times$ is an abelian group of infinite rank, we obtain that $F^\times$ is an (indecomposable) abelian group of infinite rank as well.

What can we say about the extension $F\, / \, \F_2(T)$? We claim it cannot be finite. Assume to the contrary that $d = [F\colon \F_2(T)] < \infty$. The field norm 
\[ N = N_{F/\F_2(T)} \colon F^\times \longrightarrow \F_2(T)^\times \]
is an abelian group homomorphism (see \cite[8.1.3]{Roman}) whose restriction to $\F_2(T)^\times$ is the $d$th-power map. The image of $N$, $\mathrm{im}\, N$, is a free abelian group since it is a subgroup of the free abelian group $\F_2(T)^\times$, and it is a basic fact from homological algebra that this implies $\mathrm{im}\, N$ is a summand of $F^\times$ (see \cite{DummitFoote}; we also give a direct argument in the section on valued fields below). Further, $\mathrm{im}\, N$ has infinite rank since the restriction of $N$ to $\F_2(T)^\times$ is the $d$th-power map. Thus, $F^\times$ is decomposable, a contradiction. 

This eliminates two other popular classes of fields. An {\em algebraic function field} of $n$ variables over a field $k$ is a finite extension of the field $k(x_1, \dots,x_n)$ of rational functions in $n$ variables (possibly empty set) over $k$ (this includes the {\em finitely generated fields}, where $k = \F_p$ or $\Q$). {\em Global fields} are number fields (finite extensions of $\Q$) or function fields over a finite field $\F_q$ (finite extensions of $\F_p(t)$). Such fields with characteristic 2 have decomposable multiplicative groups since they are finite extensions of fields of rational functions.

What can be said of the size of $T$? In \cite{may}, May proves that when $B$ is a finitely generated field and $A$ is an extension field generated over $B$ by algebraic elements of bounded degree, then $A^\times$ is free modulo torsion. In our situation, if $T$ is finite, then $F$ cannot be generated over $\F_2(T)$ by algebraic elements whose degrees over $\F_2(T)$ are bounded, for otherwise May's result would imply that $F^\times$ is a {\em free} indecomposable abelian group of infinite rank, which is impossible. This eliminates examples of the following sort. Let $J$ be the subfield of the algebraic closure of $\F_2(x)$ that is generated by adjoining to $\F_2(x)$ the cube root of every irreducible polynomial in $\F_2[x]$. The field $J$ is an infinite algebraic extension of $\F_2(x)$, and it must have a decomposable multiplicative group by the previous discussion.

We next discuss {\em perfect fields} in order to give an important class examples where the extension $F \, / \, \F_2(T)$ is an infinite algebraic extension by elements of unbounded degree. Every field of characteristic zero is perfect, and a field $K$ of characteristic $p > 0$ is perfect if and only if the Frobenius endomorphism $x \mapsto x^p$ is a field automorphism of $K$. We need only consider perfect fields of characteristic 2. Given a set $T$, write \[ \sqrt[2^k]{T} = \{ \sqrt[2^k]{t} \, | \, t \in T\}.\] A fundamental perfect field $F$ of characteristic 2 is a field of the form \[ \bigcup_{k \geq 0} \F_2(\sqrt[2^k]{T}) \] for some set $T$. (The field $F$ is a perfect closure of $\F_2(T)$.) Any nonzero element $Q \in F$ lies in $\F_2(\sqrt[2^k]{T})$ for some $k$, and since the squaring map is a field automorphism of $F$, \[ Q^{2^k} \in (\F_2(T))^\times \cong \bigoplus_{f \in \Delta} f^\Z.\] This implies that \[ F^\times \cong \bigoplus_{f \in \Delta} f^{\Z[1/2]}\] and hence $F$ has a decomposable multiplicative group. For perfect fields, there is a nice analogue of the fact that every field is an algebraic extension of the subfield of rational functions generated by a transcendence basis: every perfect field of characteristic 2 is an algebraic extension of a fundamental perfect field (\cite{Viet}). We do not know whether there exists such an extension with indecomposable multiplicative group.

The above discussion proves the following theorem, which summarizes what we know so far.

\begin{theorem} Let $F\, / \, \F_2$ be an infinite field and let $T$ be a transcendence basis for $F$ over $\F_2$. If $F^\times$ is indecomposable, then $T$ is nonempty and the following conditions hold. \label{newmain}
\begin{enumerate}
\item $F^\times$ is a torsion-free abelian group of infinite rank. \label{groupprops}
\item $F\, / \, \F_2(T)$ is an infinite algebraic extension.
\item If $T$ is finite, then $F$ cannot be generated over $\F_2(T)$ by algebraic elements whose degrees over $\F_2(T)$ are bounded. \label{mayitem}
\end{enumerate}
\end{theorem}

We show in \cite{cl-3} that if $G$ is an indecomposable torsion-free abelian group, then $G$ is the group of units in the group ring $\F_2[G]$. Whether $G$ is also the group of units of a field depends on whether $\F_2[G]$ has a quotient field with unit group $G$. A relatively simple example of an indecomposable torsion-free abelian group of infinite rank (and hence an example of a group that satisfies Theorem \ref{newmain} (\ref{groupprops})) is the additive group of $p$-adic integers $\Z_p$; see \cite{fuchsprob}. However, we do not know whether $\Z_p$ is the group of units of a field.

We will consider one other popular class of fields: {\em valued fields}.

\section{Valued fields} Many important classes of fields are valued fields; see \cite{vf} for a detailed introduction. Among the valued fields of characteristic 2, many can be eliminated as possible counterexamples to Conjecture \ref{mainc}.

Let $G$ be a (linearly) ordered abelian group (see \cite{vf} for the precise definition; Levi proved in \cite{Levi} that an abelian group admits such an ordering if and only if it is torsion-free). Examples include ordered subgroups of the real numbers under addition. Let $K$ be a field, and let $\map{v}{K^\times}{G}$ be a surjective group homomorphism. We may extend $v$ to a map $\map{v}{K}{G \cup \{\infty\}}$ by setting $v(0) = \infty$. The map $v$ is called a valuation on $K$ if it additionally satisfies the condition \[ v(x + y) \geq \min\{v(x), v(y)\}\] for all $x, y \in K$. (Note that $\infty \geq g$ for all $g \in G$.) The group $G$ is called the value group of $K$. A valued field with value group $\Z$ is called a discrete valued field.

Here are a few examples. The first three are variants of the same idea, but we present them separately for concreteness; they are examples of discrete valued fields.

\begin{enumerate}
\item {\bf Fields of rational functions}. Let $K$ be a field, let $X$ be a set of indeterminates, fix $x \in X$, and let $K' = K(X \setminus\{x\})$. Fix an irreducible polynomial $p(x) \in K'[x]$. For any nonzero $Q(x) \in K'(x) = K(X)$, we may write $Q(x) = p(x)^n \frac{f(x)}{g(x)}$ where $\gcd(p, f) = \gcd(p, g) = 1$. Setting $v_p(Q) = n$ determines a valuation $\map{v_p}{K(X)^\times}{\Z}$.
\item {\bf The $p$-adic valuation.} For any prime $p$, the $p$-adic valuation $v_p$ on $\Q$ is defined by $v_p(x) = n$, where $n$ is the unique integer such that $x = p^n a/b$ and $(p, a) = (p, b) = 1$. (This valuation determines an absolute value on $\Q$, and the completion of $\Q$ with respect to this absolute value produces the field of $p$-adic numbers $\Q_p$. The only nontrivial valuations on $\Q$ are the $p$-adic ones.)
\item {\bf Fields of meromorphic functions.} Let $R$ be a Riemann surface, let $M$ denote the field of meromorphic functions on $R$, and fix $p \in R$. Every nonzero $f \in M$ is locally expressible in the form $(z - p)^nh(z)$ for some holomorphic function $h$ with $h(p) \neq 0$. Setting $v(f) = n$ determines a valuation $\map{v}{M(X)^\times}{\Z}$.
\item {\bf Fields of Hahn series} (see \cite{mac}). These examples generalize rings of formal power series by allowing exponents in an arbitrary ordered abelian group $G$. (This class includes fields of formal Laurent series and fields of Puiseux series.) Let $K$ be a field. The field of Hahn series $K((G))$ consists of formal sums \[ A = \sum_{g \in G} a_g x^g \] where $a_g \in K$ and the support $\mathrm{Supp}\, A = \{g \in G \, | \, a_g \neq 0\}$ is a well-ordered set. For $A \in K((G))^\times$, $v(A) = \min\{\mathrm{Supp}\, A\}$ determines a valuation $\map{v}{K((G))^\times}{G}.$
\end{enumerate}

For a valued field $K$, a section $s$ of the valuation map $v$ is a right inverse of $v$ (i.e., a map $\map{s}{G}{K^\times}$ such that $v(s(g)) = g$ for all $g \in G$). When a section $s$ exists, the map $\map{\Phi}{G \times \ker v}{K^\times}$ defined by $\Phi(g, w) = s(g)w$ is an isomorphism (see \cite{DummitFoote} for details). This observation, in conjunction with Theorem \ref{newmain}, gives the following theorem.

\begin{theorem} Let $K$ be a valued field with nontrivial value group $G$. If the valuation on $K$ admits a section, then $K^\times \cong G \times \ker v$. In particular, if the valuation admits a section and $K^\times$ is indecomposable, then $K^\times \cong G$ and $G$ is a torsion-free abelian group of infinite rank.\end{theorem}

If the value group $G$ is a free abelian group, then a section exists: for each $g$ in a set of generators for $G$, define $s(g)$ to be an arbitrary element in $v^{-1}(\{g\})$. Since $G$ is free, this determines an abelian group homomorphism that is right inverse to $v$ by construction. Thus, any valued field whose value group is free must have a decomposable multiplicative group. This includes the discrete valued fields. For the field of Hahn series $K((G))$, the map $\map{s}{G}{K((G))^\times}$ defined by $s(g) = x^g$ is clearly a section, so the proper subgroup $x^G$ of $K((G))^\times$ is a summand, and $K((G))^\times$ is decomposable.

The unit groups of local fields have been computed (see \cite{ant}), but their decomposability may be immediately deduced from the above work without using the full unit group computation. A {\em local field} is a field that is isomorphic as a topological field to $\C$, $\R$, a finite extension of the $p$-adic numbers $\Q_p$, or a field of formal Laurent series $\F_q((X))$ over a finite field $\F_q$. We need only consider the fields of formal Laurent series of characteristic 2, which are discrete valued fields (they are fields of Hahn series with $G = \Z$) and hence have decomposable multiplicative groups.

\end{document}